%krandick.tex:  Exploring Werner Krandick's Binary Tree Jump Statistics
%%a Plain TeX file  Shalosh B. Ekhad and Doron Zeilberger ( pages)

%begin macros

\baselineskip=14pt
\parskip=10pt

\magnification=\magstephalf

\def\B{{\cal B}}

\def\1{{\overline{1}}}
\def\2{{\overline{2}}}
\parindent=0pt
\overfullrule=0in

\def\frac#1#2{{#1 \over #2}}
%\headline={\rm  \ifodd\pageno  \RightHead  \else  \LeftHead  \fi}
%\def\RightHead{\centerline{
%Title
%}}
%\def\LeftHead{ \centerline{Doron Zeilberger}}
%end macros
\centerline
{\bf 
Exploring Werner Krandick's Binary Tree Jump Statistics
}
\bigskip
\centerline
{\it Shalosh B. EKHAD and Doron ZEILBERGER}
\bigskip

{\bf Preface}

Fabrice Rouillier and  Paul Zimmermann [RZ] proposed a unified and very efficient algorithm for finding
real roots of univariate polynomials based on the good-old {\it Descartes' rule of signs}.
It improved previous algorithms due to George Collins and Alkiviadis G. Akritas, Jeremy Johnson, and Werner Krandick
(see [RZ] for references).
One reason it was so efficient was that in the process, the algorithm constructs a certain binary tree, 
that is traversed in depth-first-search and whenever there is a ``jump" it is expensive.
It turned out, that on average, there are not so many jumps, explaining the efficiently.

This motivated Werner Krandick [K] to find explicit expressions for the expectations of the statistics ``number of jumps"
and ``sum of the jump-distances" (see below for the exact definitions).
Using clever but ad hoc human reasoning, he found that they were  $\frac{n-1}{2}$ and $\frac{n(n-1)}{n+2}$ respectively (Theorem 11 of [K]).

Here we show the power of symbolic computation and experimental mathematics to do much more.
In particular, explicit expressions for the {\it variances} 
(namely $\frac{n^{2}-1}{8 n -4}$ and $\frac{2 n \left(2 n^{2}-n -1\right)}{n^{3}+7 n^{2}+16 n +12}$ respectively).
Better still, we will derive explicit expressions for the weight-enumerators of the set of full binary trees
according to these statistics (and the number of internal vertices) from which the expectation, variance,
and any number of higher moments can be easily deduced. Everything is implemented in the
Maple package {\tt Krandick.txt}, available from the front of this article:

{\tt https://sites.math.rutgers.edu/\~{}zeilberg/mamarim/mamarimhtml/krandick.html} \quad .

But first {\it definitions}.

{\bf Definitions}

$\bullet$ {\bf A (full) binary tree} is an unlabeled ordered tree, where every vertex has either $0$ children (and then it is called a {\it leaf}) or
$2$ children (and then it is called an {\it internal vertex}). A good way to define these creatures is {\it recursively}.
A binary tree has either  $0$ internal vertices (i.e. it only consists of the root), or else the root has a {\it left son} 
that is the root of a binary tree $T_L$ and a {\it right son}, that is the root of a binary tree $T_R$.

In symbols: either $T\,=\, . \, $ or $T=[T_L,T_R]$. 

$\bullet$ Let $V(T)$ be the number of internal vertices of the binary tree $T$. Note that it may be defined recursively by:
$$
V(.)=0 \quad , \quad  V([T_L,T_R])= V(T_L)+V(T_R)+1 \quad ,
$$
since by removing the root,  you lose an internal vertex.

$\bullet$ J(T) denotes the {\it number} of ``jumps" when you traverse it in {\it depth-first-search} (see [K]). It is best defined
{\it recursively} as follows
$$
J(.)=0 \quad
$$
and
$$
J([T_L,T_R])= \cases{
J(T_R), & if $T_L=.$ \quad ; \cr
J(T_L)+J(T_R)+1, \quad & otherwise. \cr}
$$

We also need an auxiliary statistic, the  {\it depth of the rightmost leaf}. It may be defined recursively as follows.

$$
D(.)=0 \quad ,
$$
and
$$
D([T_L,T_R])= 1+D(T_R) \quad .
$$

Another statistic studied in [K] is the {\it sum of jump distances}. It may be defined recursively as follows:

$$
JD(.)=0 \quad ,
$$
and
$$
JD([T_L,T_R])= JD(T_L)+JD(T_R)+D(T_L) \quad .
$$

As noticed in [K], it is readily seen that $JD(T)+D(T)=V(T)$, so in some sense, as we will see soon, the statistic $JD$ is `easier' than $J$.

{\bf Theorems}

Theorem {\bf 0}: Let $\B$ be the (infinite)  set of {\it all} binary trees, and let $f(x)$ be its {\it weight-enumerator} according to the weight $W_0$ defined by:
$W_0(T):= x^{V(T)}$. Then $f(x):=W_0(\B)$, a certain {\it formal power series}  in $x$ with {\it integer coefficients},  satisfies the quadratic equation

$$
f(x)=1 \,+ \, x\,f(x)^2 \quad .
$$

{\bf Proof}: A binary tree is either trivial, with zero internal vertices, whose weight is $x^0=1$ explaining the `$1$' on the
right side of the equation, or it has
a left tree and right tree, the $x$ in front of the second term on the right is because when you remove the root you lose
an internal vertex, and  $T_L$ and $T_R$ range all over $\B$, and since $x^{a+b}=x^a x^b$, we have it.

{\bf Comment}: Solving the quadratic equation gives
$$
f(x)= \frac{1-\sqrt{1-4x}}{2x} \quad,
$$
that thanks to Isaac Newton's binomial theorem,  implies that the number of binary trees with $n$ internal vertices, let's call it $b_n$, is given by the
good old Catalan numbers
$$
b_n=\frac{(2n)!}{n!(n+1)!}  \quad.
$$

See $[CDZ]$ for many other proofs of this result.  Also see [S] for many other combinatorial objects counted by the Catalan numbers.

{\bf Theorem  1}: 
Let $\B$ be the (infinite)  set of {\it all} binary trees.
Define a (tri-variate) weight, $W_1(T):= x^{V(T)} t^{D(T)} q^{J(T)} $,
and let $F(x,t,q):=W_1(\B)$, a certain formal power series in $x$ with coefficients that are polynomials of $t$ and $q$.
$F(x,t,q)$ satisfies the following functional equation:

$$
F(x,t,q)=1\,+ \, x\,t\,F(x,0,q)\, F(x,t,q)\,+\, xtq \left ( F(x,1,q)-F(x,0,q) \right )\, F(x,t,q)  \quad.
$$

{\bf Proof:} Follows easily from the recursive definitions of $V(T)$, $J(T)$, and $D(T)$.

{\bf  Theorem 2}: An explicit expression, in terms of `radicals', for $F(x,t,q)$ is
$$
F(x,t,q)\,=\,-\frac{-q t x +t x +\sqrt{q^{2} t^{2} x^{2}-2 q \,t^{2} x^{2}-2 q \,t^{2} x +t^{2} x^{2}-2 t^{2} x +t^{2}}+t -2}{2 \left(q t x +t^{2} x -t x -t +1\right)} \quad.
$$

{\bf Proof}: Let the right side be $G(x,t,q)$. We claim that
$$
G(x,t,q)- \left ( 1\,+ \, x\,t\,G(x,0,q)\, G(x,t,q)\,+\, xtq\cdot \left (G(x,1,q)-G(x,0,q) \right)\, G(x,t,q) \right )\, = \, 0 \quad,
$$
(check!). The theorem follows from the obvious {\bf uniqueness} of the solution of this functional equation (in the ring of formal power series in $x,t,q$).

{\bf Secret from Kitchen}: There is a sophisticated method (that we dislike!) called the {\it kernel method}, and presumably it could be done that way.
But a much better way, is to {\it hope} that in addition to the functional equation, mixing $F(x,t,q)$ and $F(x,0,q)$ and $F(x,1,q)$, it also satisfies
a {\it pure} quadratic equation with coefficients that are polynomials in $x,t$. So, using the combinatorially
derived functional equation,  we cranked out sufficiently many terms and {\it guessed} such an equation.
We found it! Of course, so far this is {\it only} a guess. Then we asked Maple to {\tt solve} it in radicals. This is still a guess. But {\bf once conjectured} it is a {\it routine verification}, that
Maple kindly did for us. See  procedure {\tt ProveJxtq()} in our Maple package {\tt Krandick.txt}.

Indeed, if you downloaded {\tt Krandick.txt} to your own laptop (that has Maple), please type:

{\tt ProveJxtq();} \quad ,

and before you know it you would get 

{\tt true}.

Note that we needed the variable $t$, corresponding to the `depth of the rightmost leaf',  in order to be able to set-up the functional equation, but we are really not interested in it!
It is only a {\it stepping stone}, a {\it catalytic variable}. At the {\it end of the day}, we can plug-in $t=1$ and get the following theorem.

{\bf Theorem  3}: 
Let $\B$ be the (infinite)  set of {\it all} binary trees.
Define a (bivariate) weight, $W_2(T):= x^{V(T)} q^{J(T)} $,
and let $H(x,q):=W_2(\B)$, a certain formal power series in $x$ with coefficients that are polynomials of $q$. We have:

$$
H(x,q)=-\frac{-q x +\sqrt{q^{2} x^{2}-2 q \,x^{2}-2 q x +x^{2}-2 x +1}+x -1}{2 q x} \quad .
$$

{\bf Theorem  4}: 
Let $\B$ be the (infinite)  set of {\it all} binary trees.
Define a (bivariate) weight, $W_3(T):= x^{V(T)} t^{D(T)} $,
and let $J(x,t):=W_3(\B)$, a certain formal power series in $x$ with coefficients that are polynomials of $t$. We have
the following functional equation:

$$
J(x,t)=1 \,+ \, x\,t\,J(x,1)J(x,t) \quad .
$$

{\bf Proof:} A member of $\B$ is either the singleton tree, `$.$', or else can be written as $T=[T_L,T_R]$.
Since $V(T)=V(T_L)+V(T_R)+1$, and $D(T)=D(T_R)+1$, the equation follows  (the left tree $T_L$ does not contribute to the $t$ part, so the variable $t$ is set to $1$).

Solving for $J(x,t)$ gives that it equals $1/(1-xtJ(x,1))$. But $J(x,1)$ is nothing but our good old $f(x)$, the generating function for the Catalan numbers. 

So we have

{\bf Theorem 5}: An explicit expression for $J(x,t)$ is
$$
J(x,t) \, = \,\frac{2}{t \sqrt{1-4 x}-t +2} \quad .
$$

{\bf Theorem  6}: 
Let $\B$ be the (infinite)  set of {\it all} binary trees.
Define a (bivariate) weight, $W_4(T):= x^{V(T)} q^{JD(T)} $,
and let $K(x,q):=W_4(\B)$, a certain formal power series in $x$ with coefficients that are polynomials of $q$. We have
the following explicit expression

$$
K(x,q) \, = \,\frac{2 q}{\sqrt{-4 q x +1}-1+2 q}  \quad .
$$

{\bf Proof}: We noticed above that $D(T)+JD(T)=V(T)$, hence $K(x,q)=J(qx,\frac{1}{q})$, and Theorem 6 follows from Theorem 5.

{\bf Moments of The Jump Statistics}

The weight-enumerator contains {\bf all} the information needed for {\it all} the moments.

In particular The generating function of the quantity 

{\it Sum of the `number of jumps'}

over all binary trees with $n$ internal vertices,  what Krandick [K] denoted by $j_n$, is the coefficient of $x^n$ in $\frac{\partial}{\partial q} H(x,q) \vert_{q=1}$ ,
that implies that the expected number of jumps is $j_n/b_n$, that happens to be  $(n-1)/2$. 

More generally, the generating function for the quantity

{\it sum of the `$r^{th}$-power of  the number of jumps'}

 over all binary trees with $n$ internal vertices,  is

the coefficient of $x^n$ in $(q\frac{\partial}{\partial q})^r H(x,q) \vert_{q=1}$ .

Calling this quantity $j_n^{(r)}$, the $r$-th moment is   $j_n^{(r)}/b_n$ \quad.

From the  usual moments, one easily derives the {\it moments about the mean}, in particular the {\it variance}.
Once we have explicit expressions for the moments about the mean (for as many as we desire), we get
the {\it scaled moments} and then we can take the limit as $n$ goes to $\infty$. To our pleasant surprise
these (at least up to the $10$-th moment) coincide with those of the normal distribution,
$0$ for odd moments, and $\frac{(2r)!}{2^r r!}$ for the $2r^{th}$ moment.
This indicates that Krandick's jump statistics is most probably {\it asymptotically normal}. Can you prove it?

{\bf Added in the new version:} Stephen Melczer and Tiadora Ruza brilliantly proved this asymptotic normality.
See their nice writeup:

{\tt https://sites.math.rutgers.edu/\~{}zeilberg/mamarim/mamarimhtml/krandickSteveTia.pdf} \quad.

It is (probably) not hard to prove that  the moments, and hence the moments about the mean, are
{\bf rational functions} of $n$,  and one can easily bound the degrees of the numerator and denominators. So why
not crank out many `data values' and fit them into rational functions? That's exactly what we did.
Notice that it is irrelevant whether we have an {\it a priori} proof that these are rational functions. Once
{\it conjectured} it is a routine (rigorous!) verification.

So we have the following experimentally derived, but fully {\it rigorizable}  theorems.

{\bf Theorem 7.1}: (first proved in [K]) The expected `number of jumps' among all binary trees
with $n$ internal vertices is
$$
\frac{n-1}{2} \quad .
$$

{\bf Theorem 7.2}:  The  variance of the `number of jumps' statistic among all binary trees
with $n$ internal vertices is
$$
\frac{n^{2}-1}{8 n -4} \quad .
$$

{\bf Theorem 7.3}:  The  kurtosis (aka `scaled fourth moment-about-the-mean') of the `number of jumps' statistic among all binary trees
with $n$ internal vertices is
$$
\frac{6 n^{3}-11 n^{2}-2 n +3}{2 n^{3}-3 n^{2}-2 n +3} \quad .
$$

(Note that it tends to $3$, as it should).

{\bf Theorem 7.4}:  The   $6^{th}$ scaled moment-about-the-mean of the `number of jumps' statistic among all binary trees
with $n$ internal vertices is
$$
\frac{60 n^{6}-300 n^{5}+391 n^{4}-20 n^{3}-82 n^{2}-16 n +15}{4 n^{6}-16 n^{5}+7 n^{4}+32 n^{3}-26 n^{2}-16 n +15} \quad ,
$$

note that it tends to $1 \cdot 3 \cdot 5=15$, as it should).

{\bf Theorem 7.5}:  The  scaled $8^{th}$ moment-about-the-mean of the `number of jumps' statistic among all binary trees
with $n$ internal vertices is
$$
\frac{840 n^{9}-7980 n^{8}+27006 n^{7}-38933 n^{6}+23070 n^{5}-6937 n^{4}+3178 n^{3}-1167 n^{2}-142 n +105}{8 n^{9}-60 n^{8}+118 n^{7}+75 n^{6}-402 n^{5}+135 n^{4}+418 n^{3}-255 n^{2}-142 n +105} \quad .
$$

(Note that it tends to $1 \cdot 3 \cdot 5 \cdot 7=105$, as it should).

For the explicit expression for the tenth scaled moment about the mean, look at the output file:

{\tt https://sites.math.rutgers.edu/\~{}zeilberg/tokhniot/oKrandick1.txt} \quad .

{\bf Moments of The Sum of Jump Distances Statistic}

This one is even more {\bf concentrated about the mean}, since as will see below, the variance tends to a constant.

Using our guessing methodology we have the following theorems.

{\bf Theorem 8.1}: (first proved in [K]) The expected `sum of jump distances'  among all binary trees
with $n$ internal vertices is
$$
\frac{n \left(n -1\right)}{n +2} \quad .
$$

{\bf Theorem 8.2}: The variance of the statistic   `sum of jump distances'  among all binary trees
with $n$ internal vertices is
$$
\frac{2 n \left(2 n^{2}-n -1\right)}{n^{3}+7 n^{2}+16 n +12} \quad,
$$

(note that it converges to $4$, hence the standard-deviation converges to $2$).

{\bf Theorem 8.3}: The {\it skewness} (aka {\it scaled third moment-about-the-mean})  of the statistic   `sum of jump distances'  among all binary trees
with $n$ internal vertices is
$$
\frac{3 \sqrt{2}\, \sqrt{\frac{\left(n^{3}-n^{2}-8 n +12\right) n}{2 n^{4}+15 n^{3}+23 n^{2}-24 n -16}}}{2} \quad,
$$

(note that it converges to $\frac{3}{2}$). In particular, this statistic is {\bf not} asymptotically normal, since for the latter to be true it should have been $0$).

{\bf Theorem 8.4}: The kurtosis of the statistic   `sum of jump distances'  among all binary trees
with $n$ internal vertices is
$$
\frac{25 n^{5}+58 n^{4}-45 n^{3}-34 n^{2}-172 n -48}{2 \left(2 n^{4}+17 n^{3}+30 n^{2}-29 n -20\right) n} \quad ,
$$

(note that it converges to $\frac{25}{4}$).

For explicit expressions for the fifth through the tenth scaled moments about the mean, look at the output file

{\tt https://sites.math.rutgers.edu/\~{}zeilberg/tokhniot/oKrandick2.txt} \quad .

{\bf Conclusion}
 
Werner Krandick used pure human cleverness to find explicit expressions for the {\it first} moments of the jump
statistics that he studied. But using  {\bf symbolic computation} and {\bf experimental mathematics},  one can go much further.
{\it We believe that this is the way to go.}

{\bf References}

[CDZ] Shaoshi Chen, Robert Dougherty-Bliss, and Doron Zeilberger, {\it $C_4$ proofs that the number of binary trees with n+1 leaves is given by the Catalan number $C_n$},
in preparation.

[K] Werner Krandick, {\it Trees and jumps and real roots},
Journal of Computational and Applied Mathematics {\bf 162} (2004), 51-55. \hfil\break
{\tt https://sites.math.rutgers.edu/\~{}zeilberg/akherim/krandick2024.pdf} \quad .

[RZ] Fabrice Rouillier and  Paul Zimmermann, {\it Efficient isolation of polynomial's real roots},
Journal of Computational and Applied Mathematics {\bf 162} (2004), 33-50. \hfil\break
{\tt https://sites.math.rutgers.edu/\~{}zeilberg/akherim/rouillier2004.pdf} \quad .

[S] Richard P. Stanley, {\it ``Catalan Numbers''}, Cambridge University Press, 2015.
\bigskip
\hrule
\bigskip

Shalosh B. Ekhad and Doron Zeilberger, Department of Mathematics, Rutgers University (New Brunswick), Hill Center-Busch Campus, 110 Frelinghuysen
Rd., Piscataway, NJ 08854-8019, USA. \hfill\break
Email: {\tt ShaloshBEkhad at gmail  dot com}, {\tt DoronZeil at gmail  dot com}   \quad .
\bigskip
{\bf Exclusively published in the Personal Journal of S.B. Ekhad and D. Zeilberger and arxiv.org}.  
\bigskip
{\bf First written: July 16, 2024.} \quad, {\bf This version: Aug. 23, 2024.} \quad,

\end